\documentclass[a4paper,12pt]{elsarticle}

\usepackage{multicol}
\usepackage{color}
\usepackage{xspace}
\usepackage{pdfwidgets}

\begin{document}

\title{Polynomial Form of the Matrix Exponential}
 \author{Daniel Gebremedhin}
 \ead{daniel1.gebremedhin@famu.edu}
 \author{Charles Weatherford}
 \ead{charles.weatherford@famu.edu}
\address{Physics Department, Florida A\&M University, Tallahassee, FL, USA.}
\date{\today}
\begin{abstract}
An algorithm  for numerically computing  the exponential of a  matrix
is presented. We  have derived a  polynomial expansion of $e^x$  by
computing it  as an  initial value  problem using  a symbolic
programming language. This algorithm  is shown  to be  comparable in
operation  count and convergence with the  state--of--the--art method
which is  based on a Pade approximation of  the exponential matrix
function. The  present polynomial form,  however, is more  reliable
because the  evaluation  requires only  linear combinations  of the
input matrix.  We also show that  the technique  used to  solve  the
differential  equation, when implemented symbolically, leads to a
rational as well as a  polynomial form of the solution function. The
rational form is the well-known diagonal Pade approximation  of
$e^x$. The polynomial form, after some rearranging to minimize
operation count, will be used to evaluate the exponential of a matrix
so as to illustrate its advantages as compared with the Pade form. 
\end{abstract}


\maketitle

\section{Introduction}
\label{sec:introduction}

The exponential function (EF) of a square matrix (matrix exponential
function--MEF) is one of the most important functions in a
computational linear algebra. It is a fundamental topic of research
pertaining to functions of matrices. Since it can be modeled as a
solution to an initial value problem, such a technique developed for
its successful computation can likely be adopted to solve other
problems of physical interest. The solution of the  time--dependent
Schr\"odinger equation (TDSE) is such a problem. \cite{tannor2007introduction}. 

Many algorithms have been developed by different authors to
numerically compute the MEF
\cite{doi:10.1137/S00361445024180}. Perhaps the most successful of
them is based on the Pade approximation (PA) of the MEF. A very
efficient implementation of the algorithm is published in
\cite{doi:10.1137/04061101X} and is adopted in programming languages
such as \emph{MATLAB} \cite{MATLAB:2014} and other numerical
scientific libraries. 

The PA of the EF is a very compact rational expansion of $e^x$ about
$x = 0$. When applied to a matrix, its evaluation must eventually
involve a matrix inversion, which can sometimes lead to a (nearly)
singular system of equations. This issue has been pointed out and
discussed in \cite{doi:10.1137/S00361445024180} in some detail. Hence,
a polynomial representation of the EF, which is as convergent and
compact as the PA can be very favorable for implementation on the
MEF. In the present work, such a polynomial form has been derived. We
will demonstrate that it approximates the MEF with similar convergence
as its PA and that it can be computed with comparable efficiency. 

In this paper, we will use variables $x$, $\tau$ etc to represent
general square matrices. References to scalar quantities should be
clear from context. We will use the term norm to mean $1$--norm
$(\|x\|_1)$ of matrix \cite{golub2013matrix}. The relative error of
matrices are also calculated based on the $1$--norm.

\section{Preliminary}
\label{sec:preliminary}

The PA can be thought of as an order $2m \rightarrow m$ numerical
economization of a polynomial representation of a function. Given the
coefficients of a general polynomial of order $2m$, it is possible to
calculate the corresponding coefficients of two other polynomials of
order $m$ so that their quotient results in a rational approximation of
the same function. This simple rearrangement usually leads to better
convergence, which makes it an interesting topic of numerical
studies. Efficient numerical algorithms are available that can
calculate these PA coefficients for an arbitrary polynomial
\cite{Press:2007:NRE:1403886}. When the same algorithm is applied to a 
Taylor series expansion of the EF, one can arrive at its diagonal
PA. In fact, these two representations of the exponential function
have known forms denoted here by $T$ and $Q$ respectively as:
\begin{equation}
  \label{eq:1}
  e^x \simeq T_{2m}(x) \quad \rightarrow \quad e^x \simeq Q_m(x) =
  \frac{R_m(x)}{R_m(-x)} 
\end{equation}
\noindent where,
\begin{equation}
  \label{eq:2}
  T_m(x) = \sum_{\mu=0}^{m}\frac{x^\mu}{\mu!}, \quad R_m(x) =
  \displaystyle\sum_{\mu=0}^m{2m - \mu \choose m}\frac{x^\mu}{\mu!}.
\end{equation}

While the advantage in reduction of operation count that $Q$ has over
$T$ in approximating the EF is immediately apparent from
Eq.~(\ref{eq:1}), simple numerical tests reveal that the resulting
increase in convergence is also quite significant. This improvement in
accuracy exhibited by the PA is generally  function dependent \cite{Press:2007:NRE:1403886}.
In the following section, we will describe an algorithm that can be
used to generate an approximate representation to the EF as a solution
of an initial value problem. By running the algorithm in symbolic
programming language we were able to observe that the method can lead
to both rational and polynomial expressions of the EF, similar to the
forms of $Q$ and $T$ above, if the appropriate set up of input parameters are
used. The resulting rational expression happens to be identical to the
diagonal PA $(Q)$ given in Eq.~(\ref{eq:1}). Under the appropriate 
circumstances, a similar algorithm can be employed to derive rational
approximations of other elementary functions also. The derived polynomial
expression however, will be the main topic of this paper, and
will be used to evaluate the MEF. In the process we will reveal
some interesting contrasts to the PA, which is the current
state--of--the--art method for evaluating the MEF. Finally numerical tests and conclusions will be
presented.
\section{Description of Algorithm}
\label{sec:descr-algor}

Let $S \equiv d/dt + p(t)$ be a first order differential operator and
$F(t)$ its solution function satisfying the following ordinary
differential equation (DE).
\begin{equation}
  \label{eq:3}
  SF(t) = q(t), \qquad t_1 \leq t \leq t_2
\end{equation}
\noindent
where, $t_1$ and $t_2$ represent end points of a particular finite
element in $t$. Let $\tau$ be a local variable with domain $-1 \leq
\tau \leq 1$ defined by the linear transformation
\begin{equation}
  \label{eq:4}
  t = \frac{1}{2} \left[ (t_2 - t_1)\tau + (t_2 + t_1) \right].
\end{equation}
\noindent
In terms of the local variable $\tau$, Eq.~(\ref{eq:3}) will be
re--written as
\begin{equation}
  \label{eq:5}
  \bar{S}\bar{F}(\tau) = \bar{q}(\tau)
\end{equation}
where the meaning of the over-bar is clear. At this point we will
expand $\bar{F}(\tau)$ in a basis set that let us explicitly fix the
initial value of the function at $t_1$ which corresponds to $\tau=-1$
in terms of the local time.
\begin{equation}
  \label{eq:6}
  \bar{F}(\tau) = \sum_{\mu = 0}^{M-1} s_\mu(\tau)B_\mu + \bar{F}(-1)
\end{equation}
where $s$ is defined in terms of Legendre polynomials of the first kind
($P$) as in \cite{Olver:2010:NHM:1830479}.
\begin{eqnarray}
  \label{eq:7}
  s_\mu(\tau) &=& \int_{-1}^\tau P_\mu(t) \rm{d}t \nonumber \\   &=&
  \frac{1}{2\mu + 1} \left[ P_{\mu+1}(\tau) - P_{\mu-1}(\tau)
  \right]
\end{eqnarray}
\noindent
with $s_0(\tau) = 1 + \tau$. The $s$-functions satisfy the following recurrence
relation.\cite{phdthesis-full-DNL}
\begin{equation}
  \label{eq:8}
  s_\mu(\tau) = \frac{1}{\mu + 1}\left[(2\mu - 1)\tau s_{\mu -
      1}(\tau) - (\mu - 2)s_{\mu - 2}(\tau)\right]
\end{equation}
\noindent
Note that $s_\mu(-1)=0$ and derivative of $s_\mu(\tau)$ is
$P_\mu(\tau)$. 

Substituting the expansion given in Eq.~(\ref{eq:6}) into
Eq.~(\ref{eq:5}), projecting from the left by $P_\nu(\tau)$ and
integrating over $\tau$ results in the following set of simultaneous
equations of size $M$: 
\begin{equation}
  \label{eq:9}
  \sum_{\mu = 0}^{M-1} \bar{\Omega}_{\nu\mu}B_{\mu} =
  \bar{\Gamma}_{\nu}
\end{equation}
\noindent
where,
\begin{eqnarray}
  \label{eq:10}
  \bar{\Omega}_{\nu\mu} &=& \int_{-1}^1 P_\nu(\tau) \bar{S}s_\mu(\tau)
  \rm{d}\tau \nonumber\\
  \bar{\Gamma}_{\nu} &=& \int_{-1}^1 P_\nu(\tau) [\bar{q}(\tau) -
  \bar{p}(\tau)\bar{F}(-1)] \rm{d}\tau.
\end{eqnarray}

From Eq.~(\ref{eq:7}) it is clear that when $p(t)$ is a constant,
$\bar{\Omega}$ is a tridiagonal matrix. After solving Eq.~(\ref{eq:9})
for $B$, we can evaluate the solution function from
Eq.~(\ref{eq:6}). The value at the end point $t_2$ is particularly
important for propagation of the solution and has the following simple
form \cite{phdthesis-full-DNL}
\begin{equation}
  \label{eq:11}
  \bar{F}(+1) = 2B_0 + \bar{F}(-1).
\end{equation}
\noindent
Note that Eq.~(\ref{eq:6}) is a polynomial of order $M$ in $\tau$ as
can be seen from the form of $s$ in Eq.~(\ref{eq:7}). 

Now, we will set up the parameters in the algorithm such that
the resulting solution function is $e^{x}$.

\subsection{Rational Approximation to the EF}
\label{sec:rati-appr-ef}

Let $F(t) = e^{tx}$, where $0 \leq t \leq 1$. The EF can be calculated
as a solution to an initial value problem by employing the above
algorithm with $p(t)=-x$, $q(t)=0$, $t_1=0$, $t_2=1$ and $F(0) =
1$. The desired expression for the EF is obtained at $F(1) = e^x$ from
Eq.~(\ref{eq:11}). By letting $x$ to be an undetermined variable, the
calculation was done in a symbolic programming language
\emph{Mathematica} \cite{Mathematica7}. For $m$ (number of basis
functions),  the resulting expression for the EF is found to be identical to
$Q_m(x)$ in Eq.~(\ref{eq:1}). It is interesting that one can calculate
the coefficients for the PA of the EF directly from the DE in this way. For the
sake of completeness, if we break the $t$ axis into, say $k$, uniform
size finite elements, the end result was found to be $F(1) =
(Q_m(x/k))^{k}$. This is consistent with the identity of the EF given
by $e^{kx} = \left( e^{x} \right)^k$. This property is more
efficiently exploited by the method of scaling and squaring, which
chooses $k$ to be a power of $2$.

\subsection{Polynomial Approximation to EF}
\label{sec:polyn-appr-ef}

Lets now drop $x$ and define $F(t) = e^t$, where $-\theta \leq t \leq
\theta$, so that we can use the above algorithm to have a polynomial
expression of the EF. $\theta$ is a positive number which determines
the domain of the result, to be chosen later based on number of basis
functions and the working machine precision. The parameters for the
algorithm will be set as: $p(t) = -1$, $q(t) = 0$, $t_1 = -\theta$,
$t_2 = \theta$ and $F(t_1) = e^{-\theta}$. The output of the algorithm
will be in a polynomial form given in Eq.~(\ref{eq:6}) which will be
valid for any $\tau$. Specifically, we will denote it by $E$ as shown
below.
\begin{equation}
  \label{eq:12}
  e^x \simeq \bar{F}(\tau = x/\theta) = E_M(x)
\end{equation}

\begin{figure}
  \includegraphics[scale=0.7]{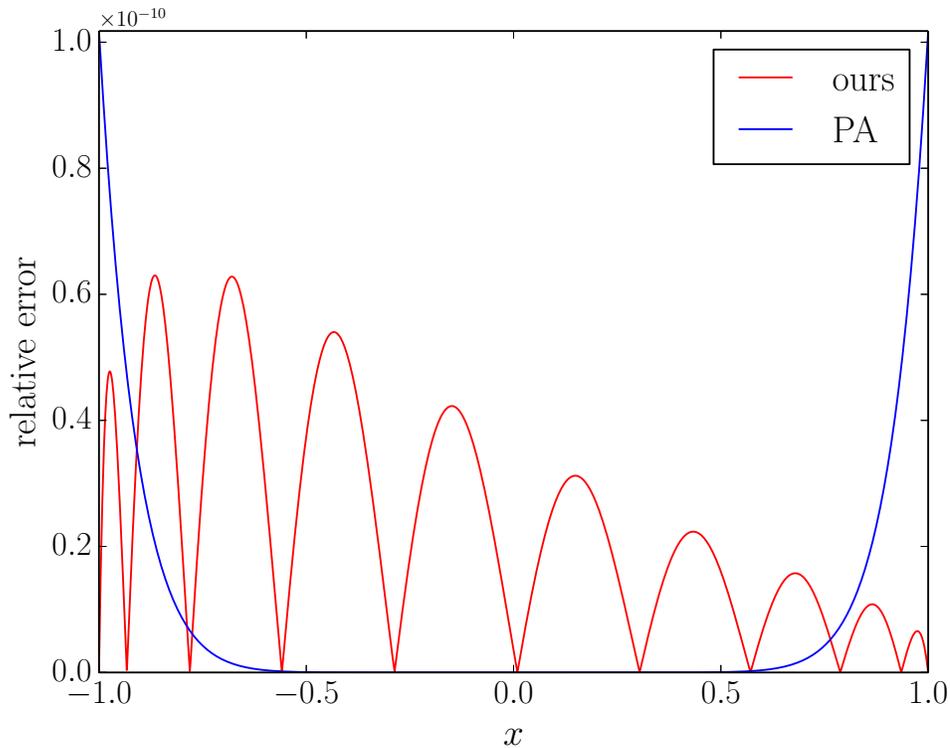}
  \caption{Relative error plots of $E_{10}(x)$ with
    $\theta = 1$, and $Q_5(x)$ in approximating $e^x$.}
  \label{fig:1}
\end{figure}

We anticipate the accuracy of $E_{2m}(x)$ to be comparable with
$Q_m(x)$ of PA. Fig.~\ref{fig:1} shows a plot of relative errors, in
approximating the scalar $e^x$, of the two methods for $m =
5$. i.e. $E_{10}(x)$ and $Q_5(x)$, where $\theta = 1$. 
The span
$\theta$ has been deliberately made too wide because lowering the
accuracy close to working precision would have altered the plots to
have the familiar random shape. Note the scale of the vertical axis.
The expansion for the PA is centered about the origin, and hence, the
U--shape; while the new polynomial is a result of a spectral method
which produces the kind of error distribution shown. Within the
domains of $\theta$s considered in this paper, the upper bound of the
relative error for the PA is generally slightly higher than the way it
is portrayed in the typical plot. This establishes the fact that
$E_M(x)$, although it is a polynomial representation, is much more
convergent than the Taylor approximation $T_M(x)$ given in
Eq.~(\ref{eq:1}). 

In passing, we notice that the choice of $p(t)=i\hat{H}$, $q(t)=0$
where, $\hat{H}$ is the Hamiltonian matrix and $i$ imaginary number,
casts the problem into a form of Schr\"odinger equation
\cite{:/content/aip/journal/jcp/128/18/10.1063/1.2916581,
  :/content/aip/journal/jcp/81/9/10.1063/1.448136, :/content/aip/journal/jcp/136/1/10.1063/1.3673320}. As will be
exhibited soon, its performance in a purely mathematical setting
promises a favorable prospect for the algorithm to be applied in
physical systems. The techniques discussed in this article can be
adopted to effectively propagate solutions to TDSE and will be
reported as soon as the details are worked out.
\section{Evaluating the Polynomial Function}
\label{sec:eval-polyn-funct}

In this section we will consider how to efficiently evaluate the EF
from its polynomial form discussed above. This can be accomplished by
minimizing the number of matrix multiplications (MMs) needed for its
evaluation. We will also choose the parameters such as the order of
expansion $M$ and the corresponding span $\theta$ so that, within the
given domain, the EF can be calculated close to the working machine
precision. 

\subsection{Product Form}
\label{sec:product-form}

Constructing $E_M(x)$ using the recursion relation shown in
Eq.~(\ref{eq:8}) can only be done by $(M - 1)$ matrix
multiplications. We need to lower the number of MMs by at least a
factor of $4$ in order to be competitive with present methods for
similar accuracy. In \cite{doi:10.1137/04061101X}, the MMs were lowered
primarily because both of the polynomials in the quotient of the PA
are already of order $M/2$. In this paper, we will
rearrange our polynomial, not as quotient, but as a product of other
polynomials of smaller order. This can be efficiently done by
calculating all the roots $x_1, x_2, \ldots x_M$ of the polynomial
given in Eq.~(\ref{eq:6}), after substituting for $\tau =
x/\theta$. This allows us to rewrite it as a product of $(x - x_1)
(x - x_2) \ldots (x - x_M)$. Generally, roots of a polynomial can be
complex numbers, in which case, they always are complex conjugate
pairs. So we need to inflate those terms that belong to conjugate
pairs, say, $x_1^*$ and $x_1$, into real quadratic expressions as $(x
- x_1)(x - x_1^*) = \left[ x^2 - 2(x_1 + x_1^*) x + x_1^*x_1
\right]$. Hence, using this procedure, we can generally write the EF
in the form
\begin{equation}
  \label{eq:13}
  E_{M}(x) =  \alpha \prod_{i = 1}^{M/2} \left( \displaystyle
    \sum_{j=0}^2 c_{ij} x^j \right) 
\end{equation}
\noindent where all the coefficients $c_{ij}$ are now real
numbers. $M$ is assumed to be an even integer. $\alpha$ takes the
value of the leading coefficient in Eq.~(\ref{eq:6}), i.e., the
coefficient of $x^M$. Now, by storing $x^2$ we can construct $E$ by a
total of $M/2$ multiplications which is clearly an improvement over $M
- 1$. We can further consider polynomials of order $4$,
$6$ etc., and seek the most economical arrangement. Finding the
order of the polynomial in the above product that requires the least
number of MMs is essentially an optimization problem of a general case
given by,
\begin{equation}
  \label{eq:14}
  E_{M}(x) =  \alpha \prod_{i = 1}^{m'} \left( \displaystyle
    \sum_{j=0}^{m} c_{ij} x^j \right) 
\end{equation}
\noindent where, $M = m'm$. Storing $x^2, \ldots , x^{m}$, which takes
$(m - 1)$ multiplies, allows us to construct any of the polynomials in
the above product. Then another $(m' - 1)$ multiplications are needed
in order to complete the evaluation of $E$. Hence, the total number of
MMs required is given by the function $(M/m + m - 2)$. The minimum of
this function occurs at $m = m' = \sqrt{M}$, giving a total of
$2(\sqrt{M} - 1)$ multiplications. Apparently, convenient values for
the order of the polynomial are squares of even integers $M = 2^2,
4^2, 6^2, \ldots, (2m)^2$, respectively requiring $2, 6, 10, \ldots,
2(2m - 1)$ multiplications. In this particular choice of $M$, a unit
increase in $m$ always raises the number of MMs by $4$.

Note that rearranging a polynomial into a product of other polynomials
of lower order, as discussed above, leads to an exactly equivalent
expression unlike the rational form of PA which can generally alter
(usually for the better) the convergence of the corresponding
polynomial form.

Once the orders of the polynomial $M$ to be used have been selected, we
need to fix the corresponding  $\theta$ so that we can calculate the
required coefficients $c$.

\subsection{Choice of Span $\theta$}
\label{sec:choice-span-theta}

In a recent paper \cite{PhysRevE.89.053319}, we have defined an
adaptive finite element step size choice, which is based on Taylor
series expansion, that is effective for solving differential
equations. When the method is applied to the EF the result is a step
size of $1.5$. We will use this quantity as a unit of measure of span and
choose from values of $\theta = 0.75m, \, m = 1, 2, \ldots$ By making
relative error plots of the scalar function similar to
Fig.~\ref{fig:1}, the largest $\theta$ value with error bounds
reasonably within the required machine precision has been selected.

Table~\ref{tab:table1} shows a summary of $M$ and $\theta$ values
considered in the present work.
%
%
\begin{table}
\caption{Values of parameters for different orders $M$ are shown. The
    $2^{nd}$ column is the number of multiplications required with
    comparison from $\pi_{M/2}$ of
    \cite{doi:10.1137/04061101X}. $\theta$ is the value of span for
    the corresponding machine precision $\epsilon$. Equivalent values
    in \cite{doi:10.1137/04061101X} are displayed as $\theta_{M/2}$}
\begin{center}\footnotesize
\renewcommand{\arraystretch}{1.3}
\begin{tabular}{|c|c|c|c|c|c|}\hline
      \textrm{$M$}& \textrm{$2(\sqrt{M} - 1)$}&
      \textrm{$\theta$}&
      \textrm{$\epsilon$} &
      \textrm{$\pi_{M/2}$} &
      \textrm{$\theta_{M/2}$} \\ \hline
      $16$ & 6 & 1.5 & $2^{-52}$ & 5 & 1.5 \\
      $36$ & 10 & 9.75 & $2^{-52}$ & 8 & 11.0 \\ 
      $64$ & 14 & 20.25 & $2^{-52}$ & -- & -- \\ 
      $64$ & 14 & 12 & $2^{-112}$ & -- & -- \\ \hline
\end{tabular}
\end{center}
\label{tab:table1}
\end{table}%

We have also included a span for $M = 64$ suitable for calculations in
quadruple precision.
Its output will be used as an exact value of the EF for
comparison purposes. Note that the method of  \cite{doi:10.1137/04061101X} still has
to solve the resulting matrix equation after constructing the matrices
in the quotient of PA by performing the indicated $\pi_{M/2}$ MMs.

For a given order $M$, all the parameters $[\alpha, c_{ij}, \theta]$
need to be calculated only once and stored. All of these
calculations that are necessary to determine the final set of parameters
have been done using exact symbolic calculations in
\emph{Mathematica}. In Table~\ref{tab:table2} the first $19$ figures
of such a result for $M = 16$ are displayed for demonstration purposes.
\begin{table}[htbp]
  \caption{\label{tab:table2}%
    Numerical values of the indicated parameters for $M$ = 16. The
    first $19$ figures are shown. All the leading coefficients are
    unity $(c_{i4} = 1)$. The number in square bracket signifies power of
    $10$.
    }
\begin{center}\footnotesize
\renewcommand{\arraystretch}{1.3}
\begin{tabular}{|c|c|}\hline
      parameter & value \\ \hline
      $\alpha$ & 4.955887515892002289[-14] \\
      $c_{13}$ & -4.881331340410683266 \\
      $c_{12}$ & -14.86233950714664427 \\
      $c_{11}$ & 862.0738730089864644 \\
      $c_{10}$ & 3599.994262347704951 \\
      $c_{23}$ & 7.763092503482958289 \\
      $c_{22}$ & 77.58934041908401266 \\
      $c_{21}$ & 430.8068649851425321 \\
      $c_{20}$ & 1693.461215815646064 \\
      $c_{33}$ & 9.794888991082968084 \\
      $c_{32}$ & 98.78409444643527097 \\
      $c_{31}$ & 387.7896702475912482 \\
      $c_{30}$ & 1478.920917621023984 \\
      $c_{43}$ & 3.323349845844756893 \\
      $c_{42}$ & 37.31797993128430013 \\
      $c_{41}$ & 545.9089563171489062 \\
      $c_{40}$ & 2237.981769593417334 \\
      $\theta$ & 1.5 \\ \hline
\end{tabular}
\end{center}
\end{table}
\section{Numerical Tests}
\label{sec:numerical-tests}

We have extensively tested the results of $E_M(x)$ on different kinds
of matrices. We will present the results of three sets of examples
herein. In all test cases we will show a comparison with the \emph{expm}
function of \emph{MATLAB} which clearly asserts that it implements the PA
according to \cite{doi:10.1137/04061101X}. Following \emph{expm}, no
preconditioning of the input matrices such as trace reduction or
matrix balancing suggested in \cite{doi:10.1137/04061101X} has been
done for comparison purpose. 

Matrices with norms less than $1.5$ and $9.75$ will be handled by $M =
16$ and $36$ respectively, while the ones with higher norms will be
calculated by $M = 36$ with proper use of scaling and squaring to
lower the norm of the matrices to below $\theta = 9.75$. As
mentioned earlier, relative errors will be calculated with reference to
the output of $E_M(x)$ in quadruple precision using parameters shown in
the last row of Table~\ref{tab:table1}.

When elements of the resulting EF of the matrix overflow beyond
$2^{1024}$ in absolute value, the corresponding relative errors are
shown as $1$ in the plots. Similarly, relative errors less than
$10^{-17}$ are overwritten to $10^{-17}$. The programming has been
done in Fortran 95, using GNU gcc version 4.9.1 compiler, in a 2.5 GHz
Intel Core i7 MacBook Pro laptop computer.

\subsubsection{The Matrix Computation Toolbox}
\label{sec:matr-comp-toolb}

\begin{figure}
  \includegraphics[scale=0.7]{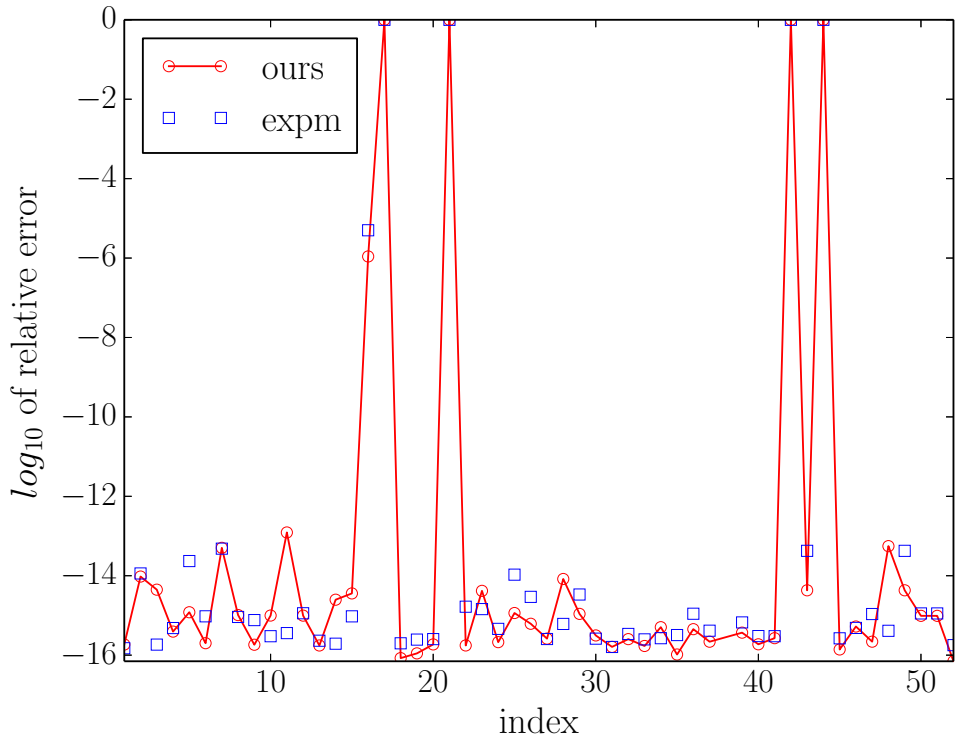}
  \caption{Test results for $51$ real matrices of size
    $8$ constructed using the matrix computation toolbox. $Log_{10}$
    of the relative errors are shown.}
  \label{fig:2}
\end{figure}

The first test matrices are constructed using the subroutines given in
the matrix computation toolbox
\footnote{http://www.maths.manchester.ac.uk/~higham/mctoolbox/} and
\cite{Higham:2002:ASN}. The size of all matrices has been set to be $8
\times 8$. Fig.~\ref{fig:2} shows the result of the numerical test
with indexes adopted from their catalog. Matrices $17, 21, 42, 44$ are
beyond the scope of the test due to overflow. On the rest of the
matrices, there is a striking qualitative similarity in accuracy
between the two functions. 

\subsubsection{Matrix Market}
\label{sec:matrix-market}

\begin{figure}
  \includegraphics[scale=0.7]{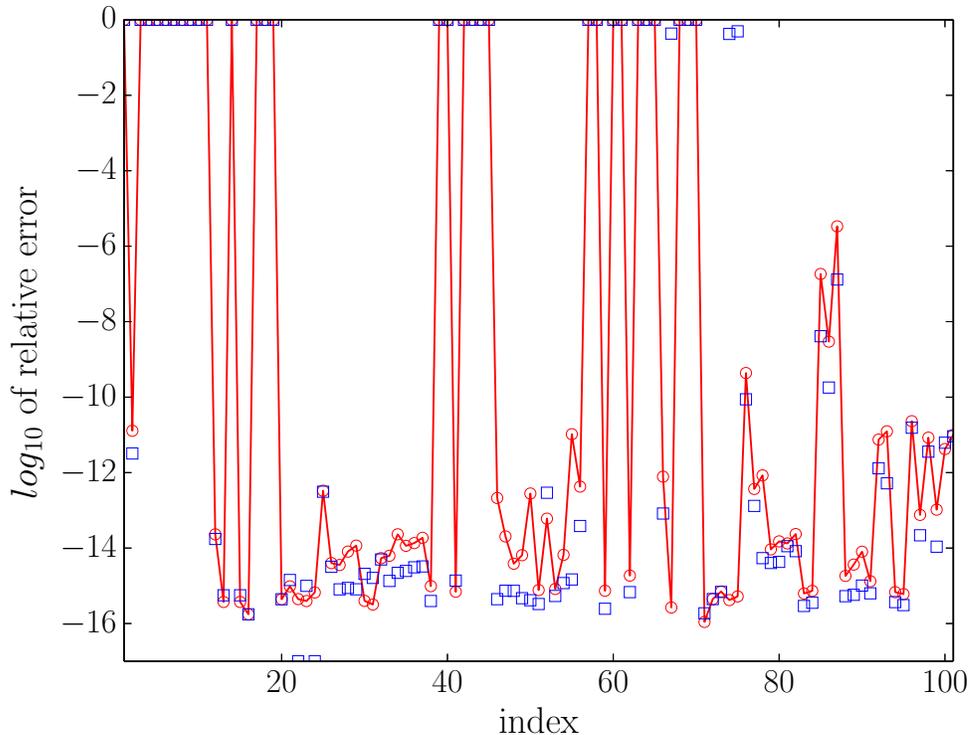}
  \caption{Test results for $101$ real matrices of
    sizes up to $500$, downloaded from matrix market. $Log_{10}$
    of the relative errors are shown.} 
  \label{fig:3}
\end{figure}

The second set of test matrices have been taken out of the matrix
market in \footnote{http://math.nist.gov/MatrixMarket/index.html} and
\cite{Boisvert:1997:MMW:265834.265854}. A query has been submitted for
real and square matrices of size up to $500$. This returned a set of
$101$ different kinds of matrices all of which has been tested
here. Relative errors of \emph{expm} and ours is shown in
Fig.~\ref{fig:3}. In both functions, $30$ matrices overflow upon evaluation. 
The qualitative similarity in accuracy between the two
methods is consistent here as well, except for three matrices labeled
with indexes $67, 74, 75$, in which there is a clear difference in
accuracy in favor of our method. Those matrices are respectively named
\lq mhd416b\rq, \lq plat362\rq \, \& \lq plskz362\rq \, in the matrix
market and have sizes $416, 362$ \& $362$. These three matrices
commonly have low norms and high sparsity, which among other reasons,
leads us to suspect that they might have posed an ill--conditioned
matrix during the matrix inversion step of the PA. To see if sparsity
is an issue, we have made a test on collection of matrices
which are mainly sparse.  

\subsubsection{UF Sparse Matrix Collection}
\label{sec:uf-sparse-matrix}

\begin{figure}
  \includegraphics[scale=0.7]{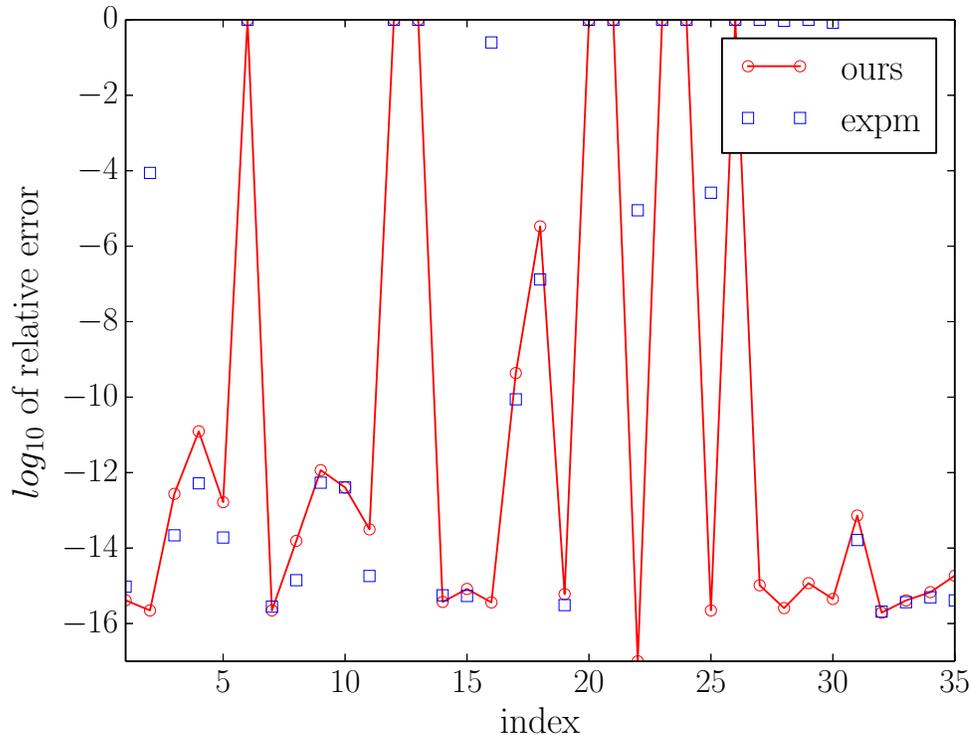}
  \caption{Test results for $35$ real matrices of
    sizes up to $100$, downloaded from UF sparse matrix
    collection. $Log_{10}$ of the relative errors are shown.} 
  \label{fig:4}
\end{figure}

The last test is downloaded from University of Florida sparse matrix
collection
\footnote{http://www.cise.ufl.edu/research/sparse/matrices/}
\cite{Davis:2011:UFS}. We downloaded and tested all real square
matrices with sizes up to $100$. The number of such matrices was 
$35$. The result of the test is plotted in Fig.~\ref{fig:4}. 
The plot clearly exposes the weakness of the \emph{expm} function in
addressing sparse matrices. The matrix labeled by index $27$, named
\lq dbGD97\_b\rq \, in the collection, for example, completely blows
up to \emph{Inf} and/or \emph{NaN}, when evaluated by \emph{expm}
which can only be explained by a singular matrix during an LU
decomposition process.  

It is not clear how to {\it a priori} identify what kind of input matrix
will eventually lead to a poorly conditioned matrix which will compromise the 
inversion step in the PA. This makes our method more reliable 
because it merely involves linear combinations of the input
matrix. Note that the above tests are exhaustive in the sense that all
the resulting matrices that fulfill the mentioned search criteria are
considered, and there are no cases where \emph{expm}
outperforms ours other than what is shown in the plots.  

Finally, although not shown here, by changing the horizontal axis in
the relative error plots to be the norm of the input matrix (instead of index), we were able to see that there is a compelling
correlation between the two. Even with the powerful method of scaling
and squaring in place, very high norm of the input matrix is known to
be a challenge inherent to computing the EF, which seems to be the
case with our algorithm as well.

\section{Conclusion}
\label{sec:conclusion}

We have derived a numerical algorithm that calculates the EF of
matrix. The EF is given in a polynomial form, which we have shown how
it can be evaluated by a minimal number of MMs. Sparsity of matrices
can be exploited element--wise during evaluation of MMs, which makes
our method more so efficient. Matrices with very high norms correspond
to longer propagation of solution, which naturally compounds error
growth as long as we are working with finite precision. It is those
matrices with very high norms that were challenging to our algorithm.  

Generally, in this work we have implemented a numerical method that
enabled the derivation of both rational and polynomial expansions of
an important mathematical function from its DE. This algorithm can
readily be adopted to a more realistic dynamical systems that can be
modeled as an evolution of initial value problem
\cite{PhysRevE.89.053319}.  

The similarity of results shown on the test plots and other aspect of
the calculations indicate that our polynomial form is indeed
complementary to PA. It is informative to see that they are both
solutions to the same DE attainable via a simple technique. But the
polynomial form is based solely on MM and avoids the matrix inversion
process altogether, and hence, as the plotted test results indicate,
ours is the less dubious one. 

\section{Acknowledgement}
\label{sec:acknowledgement}
DHG and CAW were partially supported by the Department of Energy, National Nuclear Security Administration, under Award Number(s) DE-NA0002630. CAW was also supported in part by the Defense Threat Reduction Agency.

\section{REFERENCES}
\label{sec:references}
\bibliographystyle{elsarticle-num-names}
\bibliography{JCPmatexp}

\end{document}